\pdfoutput=1
\RequirePackage{ifpdf}
\ifpdf 
\documentclass[pdftex]{sigma}
\else
\documentclass{sigma}
\fi

\numberwithin{equation}{section}

\newtheorem{Theorem}{Theorem}[section]
\newtheorem{Corollary}[Theorem]{Corollary}
\newtheorem{Proposition}[Theorem]{Proposition}
 { \theoremstyle{definition}
\newtheorem{Definition}[Theorem]{Definition}
\newtheorem{Remark}[Theorem]{Remark} }

\begin{document}

\newcommand{\arXivNumber}{1412.4721}

\allowdisplaybreaks

\renewcommand{\PaperNumber}{027}

\FirstPageHeading

\ShortArticleName{An Integrability Condition for Simple Lie Groups~II}

\ArticleName{An Integrability Condition for Simple Lie Groups~II}

\Author{Maung MIN-OO}
\AuthorNameForHeading{M.~Min-Oo}
\Address{Department of Mathematics \& Statistics,  McMaster University, Hamilton, Canada}
\Email{\href{mailto:minoo@mcmaster.ca}{minoo@mcmaster.ca}}

\ArticleDates{Received December 17, 2014, in f\/inal form March 26, 2015; Published online April 01, 2015}

\Abstract{It is shown that a simple Lie group $G$ ($ \neq {\rm SL}_2$) can be locally characterised by an integrability condition on an $\operatorname{Aut}(\mathfrak{g})$ structure on the tangent bundle, where $\operatorname{Aut}(\mathfrak{g})$ is the automorphism group of the Lie algebra of~$G$. The integrability condition is the vanishing of a torsion tensor of type~$(1,2)$. This is a slight improvement of an earlier result proved in
[Min-Oo M., Ruh E.A., in
Dif\/ferential Geometry and Complex Analysis, Springer, Berlin, 1985, 205--211].}


\Classification{53C10; 53C30}

\section{Introduction}

This is a very short addendum to a paper that I wrote with E.A.~Ruh in 1985 \cite{mr}, where we characterized a simple Lie group $G$ ($ \neq {\rm SL}_2$) by an integrability condition on an $ \operatorname{Aut}(\mathfrak{g})$ structure on the tangent bundle, where $\operatorname{Aut}(\mathfrak{g})$ is the automorphism group of the Lie algebra of~$G$. The condition was on a tensor of type $(1,3)$. In this paper we derive a slightly improved version of that result by using a more natural integrability condition involving a tensor of type $(1,2)$, namely the ``canonical'' torsion for the given structure.

The classical theorems of Sophus Lie characterise what we now call Lie groups (locally) by their Lie algebras, which are vector f\/ields whose Lie brackets are ``constant'' (structure constants).
Nowadays, geometric structures are usually described by inf\/initesimal tensorial objects such as a metric or an almost-complex structure def\/ining what is known as a $G$-structure (a reduction of the frame bundle to a subgroup). Our result gives a characterisation of a simple (local) Lie group $ \neq {\rm SL}_2$) by the torsion of the appropriate $G$-structure, namely the automorphism group of its Lie algebra, which in this case is the same as~$G$, up to coverings and connected components. Our result is dif\/ferent from the classical theorems because it does not hold for a~general Lie group, only for simple Lie groups of rank~$\geq 2$. The main reason is that we use the deep result about holonomy groups by Berger~\cite{b} and Simons~\cite{s}. More precisely, we need the purely algebraic version proved by Simons~\cite{s} about the rigidity of holonomy groups of rank $\geq 2$, acting non-transitively on the unit sphere (that is the reason why we have to exclude ${\rm SL}_2$, which is of rank one). One can f\/ind an independent simple algebraic proof of Simon's result for simple Lie groups of rank $\geq 2$ in our earlier paper~\cite{mr}.

\section{The result}

Let $M$ an $n$-dimensional manifold and $\mathfrak{g}$ a Lie algebra of the same dimension. An $\operatorname{Aut}(\mathfrak{g})$-structure, in the sense of $G$-structures is a reduction of the frame bundle of $M$ to the subgroup $\operatorname{Aut}(\mathfrak{g}) \subset  {\rm GL}(n;{\mathbb R})$. Such a structure is determined by a skew-symmetric tensor $T: TM \otimes TM \longrightarrow TM $ of type $(1,2)$ satisfying the Jacobi identity:
\begin{gather}
T(X,Y) + T(Y,X) = 0,\nonumber\\
T(X,T(Y,Z)) + T(Y,T(Z,X)) + T(Z,T(X,Y)) = 0.\label{eq2.2}
\end{gather}
$T$ def\/ines a Lie algebra structure for the tangent space at each point.

The $\operatorname{Aut}({\mathfrak g})$-structure is given by all frames $ u: {\mathfrak g} \cong T_p M$ satisfying
\begin{gather*}
 u([A , B]) = T(u(A) , u(B))
 \end{gather*}
and $G$ acts through inner automorphisms:  $u   g  = u \circ {\rm ad}_g$.
If $u$ is a moving frame (a local section of the principal bundle), then we do not assume that the Lie bracket of vector f\/ields
$ [u(A), u(B)] $ is the same as $ u([A, B])$ (otherwise we will be just reproving Lie's original theorems).

For a semi-simple Lie algebra ${\mathfrak g}$, the f\/irst Lie algebra cohomology group $H^1({\mathfrak g} ; {\rm ad})$
with values in the adjoint representation vanishes and so almost all automorphisms are inner automorphisms, in the sense that  the inner automorphisms form a normal subgroup of f\/inite index. In fact, $H^2({\mathfrak g} ; {\rm ad})$ also vanishes
(see~\cite{j}) and this is crucial for our proof. A semi-simple real Lie algebra~${\mathfrak g}$ has a natural bi-invariant non-degenerate inner product (not necessarily positive-def\/inite) given by the Killing form $\langle X , Y \rangle = - B_{{\mathfrak g}}(X,Y)$.
This def\/ines a (pseudo-)Riemannian metric on $M$, which we shall again denote by $\langle\;,\;   \rangle$.
Let $D$ denote the Levi-Civita connection of this metric. In~\cite{mr}, we looked at the following tensor $d^D T$ of type $(1,3)$
\begin{gather*}
d^D T( X,Y,Z) = (D_X T) (Y , Z) + (D_Y T) (Z , X) + (D_Z T) (X , Y)
\end{gather*}
  and we proved the following result
  \begin{Theorem}[\cite{mr}]\label{theorem2.1}
  Let ${\mathfrak g}$ be a semi-simple real Lie algebra different from ${\mathfrak{sl}}_2({\mathbb C})$ or any of its real forms. Then $d^DT = 0$ iff $(M, \langle \;, \;\rangle)$ is either locally isometric to a Lie group~$G$
  $($with Lie algebra ${\mathfrak g})$ with its unique $($up to a constant factor$)$ bi-invariant metric, or else is flat.
  \end{Theorem}

  The purpose of this paper is to introduce a more natural integrability condition. First we need to def\/ine ``dual'' bases
  $\{ e_k\}$ and $\{ e^j\}$ for ${\mathfrak g}$ satisfying $\langle e_k ,  e^j \rangle = \delta _k^j$ (this is standard procedure in this theory):

 Case 1. If ${\mathfrak g}$ is compact, then the metric $\langle \;,\; \rangle $ is positive def\/inite and we just set $e_k = e^k$, where
  $\{ e_k\}$ is an orthonormal base for the metric.

  Case 2. If ${\mathfrak g}$ is not compact, we split ${\mathfrak g} = {\mathfrak k} \oplus {\mathfrak m}$ where the metric $\langle\;  ,\; \rangle$ is positive-def\/inite on ${\mathfrak k}$ (the maximal compact subalgebra)
  and negative-def\/inite on the vector space ${\mathfrak m}$. On ${\mathfrak k}$ we def\/ine $e_k = e^k$, where $\{ e_k \}$ is an orthonormal base. On ${\mathfrak m}$, we set $e^j = - e_j$.

These def\/ine (local) frames for the $\operatorname{Aut}({\mathfrak g})$-structure on $M$.

We now def\/ine our integrability tensor $\tau$ as follows
\begin{gather*}
 \tau(X,Y) = \sum_k   \big( T\big( e_k, (D_X T) \big(Y, e^k\big) \big)  - T\big( e_k, (D_Y T ) \big(X,e^k\big) \big) \big).
\end{gather*}

This might look a bit complicated at f\/irst sight, but it is a tensor of type $(1,2)$, so it is more in line with many other integrability conditions.
It is clear that $\tau$ is a well def\/ined tensor on~$M$ (independent of the choice of the dual base), given the $\operatorname{Aut}({\mathfrak g})$-structure.
The main reason why~$\tau$ is the correct tensor is because it is the torsion of the following connection:

\begin{Definition}\label{definition2.2}
\begin{gather*}
 \nabla _X Y = D_X Y + A_X Y,
 \end{gather*}
where
\begin{gather*}
  A_X Y = \sum_k T\big( e_k, (D_X T) \big(Y, e^k\big) \big).
  \end{gather*}
\end{Definition}

Now here is the key fact about this connection:

\begin{Proposition}\label{proposition2.3}
$T$ is parallel with respect to the connection $\nabla$, i.e.,
\begin{gather*}
\nabla T = 0.
\end{gather*}
\end{Proposition}

Since the (pseudo-)Riemannian metric $\langle \; ,  \; \rangle$ on $M$ is def\/ined algebraically using only the Killing form of the Lie algebra structure on each tangent space $T_pM$, given by the tensor $T$, we have the following corollaries:

\begin{Corollary}\label{corollary2.4}
$\nabla$ is a metric connection with torsion $\tau$:
\begin{gather*}
\langle A_X Y , Z  \rangle +  \langle A_X Z , Y  \rangle = 0.
\end{gather*}
\end{Corollary}

\begin{Corollary}\label{corollary2.5}
\begin{gather*}
\tau =0 \quad \Leftrightarrow \quad A = 0 \quad \Leftrightarrow \quad DT = 0 \quad \Leftrightarrow \quad d^D T = 0.
\end{gather*}
\end{Corollary}

\begin{Remark}\label{remark2.6}
 In the language of \cite{ss} or~\cite{m}, $\tau$ is called the torsion of the $\operatorname{Aut}({\mathfrak g})$-structure, which is an invariant that depends only on the $\operatorname{Aut}({\mathfrak g})$-structure and $\nabla$ is the
 explicit connection that realises this torsion.
 \end{Remark}

So in view of Theorem~\ref{theorem2.1}, proved in the earlier paper \cite{mr}, Proposition~\ref{proposition2.3} now implies the main result of this short note:

\begin{Theorem}\label{theorem2.7}
  Let ${\mathfrak g}$ be a semi-simple real Lie algebra different from ${\mathfrak{sl}}_2({\mathbb C})$ or any of its real forms. Then $\tau = 0$ iff $(M, \langle \;, \;\rangle)$ is either locally isometric to a Lie group $G$
  $($with Lie algebra ${\mathfrak g})$ with its unique $($up to a constant factor$)$ bi-invariant metric, or else is flat.
  \end{Theorem}

\begin{proof}[Proof of Proposition~\ref{proposition2.3}]
Let $ v \in T_p M$, and let us denote for simplicity the the covariant derivative $D_v T$ by $T^{\prime}$. Then by dif\/ferentiating the Jacobi identity~(\ref{eq2.2}), we obtain
\begin{gather*}
 T(   T^{\prime}(X,Y) ,   Z   ) + T(   T^{\prime}(Y,Z) ,  X  ) +  T(   T^{\prime}(Z,X) ,   Y   )\\
\qquad{} +  T^{\prime}   (T(X,Y) ,   Z  )  +  T^{\prime}   (T(Y,Z) ,   X  )  + T^{\prime}   (T(Z,X) ,   Y  ) = 0.
\end{gather*}

In the language Lie algebra cohomology, this is simply a cocycle condition on $T^{\prime}$ can simply be expresses as
\begin{gather*}
d_{\rm ad}  T ^{\prime}  = 0,
\end{gather*}
where $d_{\rm ad}$ is the co-boundary operator for two-forms with values in the adjoint representation. It is well known \cite{j} that $H^2( {\mathfrak g} ; {\rm ad}) \cong {0}$. Hence
\begin{gather*}
d_{\rm ad}   T ^{\prime}  = 0  \Rightarrow
T ^{\prime}  =  d_{\rm ad}  A
\end{gather*}
for some $A$,  where $A$ is a ${\mathfrak g}$-valued one-form. i.e.\ a ``gauge transformation'' $A: T_p M \rightarrow {\mathfrak g} \cong T_p M$.
In fact it can be easily checked that $A$ can be explicitly written as (see \cite[p.~90]{j}):
\begin{gather*}
A(X) = \sum_k  T\big( e_k, T ^{\prime} \big(X,e^k\big)  \big).
\end{gather*}
This is in fact the ``optimal'' choice of $A$ given the metric structure, in the sense of Hodge theory, since $A$ satisf\/ies $d_{\rm ad}^{\star}A = 0$, where
$d_{\rm ad}^{\star}$ is the adjoint of $d_{\rm ad}$
with respect to the natural bi-invariant metric. We have then
\begin{gather*}
T ^{\prime}   (X,Y)  =  (d_{\rm ad}   A)(X,Y) = T(AX , Y) +T(X , AY) - A(T(X,Y)),
\end{gather*}
so
\begin{gather*}
(\nabla_v T)(X,Y)  =  (D_v T)(X,Y) + A_v( T(X,Y)) - T(A_v X, Y) - T(X, A_v Y) \\
\hphantom{(\nabla_v T)(X,Y)} {} =  T^{\prime}   (X,Y)  -  d_{\rm ad}   A(X,Y) = 0.\tag*{\qed}
\end{gather*}
  \renewcommand{\qed}{}
\end{proof}

\pdfbookmark[1]{References}{ref}
\LastPageEnding

\end{document}